\documentclass[12pt,a4paper,draft]{amsart}
\usepackage{color}
\usepackage[colorlinks=true]{hyperref}
\hypersetup{urlcolor=blue, citecolor=blue, draft=false}

\theoremstyle{remark}

\theoremstyle{remark}\newtheorem{remark}{Remark}[section]

\begin{document}

\begin{center}
{\bf IMPROVED CRITERIA FOR GEOMETRIC PROPERTIES OF CONFLUENT HYPERGEOMETRIC FUNCTIONS}\\

\bigskip
Jugal Kishore PRAJAPAT \qquad Naveen KUMARI\\
Department of Mathematics, Central University of Rajasthan, \\
Bandarsindri-305817, Kishangarh, Rajasthan, India\\

Email: jkprajapat@curaj.ac.in, naveenpuniya9@gmail.com
\end{center}

\vspace{1cm}
\noindent
{\bf Abstract:} In the paper, we establish sufficient conditions for several geometric properties of confluent hypergeometric functions. Using differential subordination techniques, we derive parameter conditions ensuring lower bounds for the real part of the confluent hypergeometric function and sectorial subordination of its derivative. We also obtain criteria for strong convexity and strong starlikeness of prescribed order. In addition, we investigate the close-to-convexity and starlikeness of the confluent hypergeometric functions. Particular cases yield corresponding results for generalized Laguerre polynomials. The obtained criteria are compared with existing sufficient conditions, and their implications for the geometric theory of special functions are discussed.

\medskip
\noindent
{\bf Keywords:} Confluent hypergeometric function; differential subordination; strong starlikeness; strong convexity; close-to-convexity; generalized Laguerre polynomials.

\medskip
\noindent
{\bf Mathematics Subject Classification (2020):} 33C15, 30C45, 30C80.

\section{Introduction}
The confluent hypergeometric function of the first kind, also called Kummer's function, is defined by
\begin{equation}\label{chf-1}
\Phi(a;c;z)={}_1F_1(a;c;z)=\sum_{k=0}^{\infty} \frac{(a)_k}{(c)_k}\frac{z^k}{k!},
\end{equation}
where $a,c\in\mathbb C$, $c\notin{0,-1,-2,\ldots}$, and $(a)_n$ denotes the Pochhammer symbol. The series converges for every $z\in\mathbb C$. The function $\Phi(a;c;z)$ satisfies Kummer's differential equation
\begin{equation}\label{chf-2}
z \, \Phi''(z)+(c-z)\Phi'(z)-a\Phi(z)=0
\end{equation}
Several classical special functions admit representations in terms of $\Phi$, including Laguerre polynomials, Whittaker functions, and modified Bessel functions.

The geometric properties of confluent hypergeometric functions have been investigated using differential subordination and coefficient techniques. Merkes and Scott \cite{Merkes} studied starlikeness properties of hypergeometric functions, while Miller and Mocanu \cite{miller2} derived sufficient conditions for confluent hypergeometric functions to have positive real part, starlike, and convex in $\Delta$. Ponnusamy and Vuorinen \cite{Ponnu} obtained parameter conditions for univalence, convexity, and close-to-convexity. Bohra and Ravichandran \cite{nbr} investigated strong starlikeness and strong convexity. Confluent hypergeometric functions also arise in problems from mathematical physics, engineering, and applied sciences; see \cite{Arzikulov, Berisha, Gil, Hasanov, Lin, Tudela}. Further results on their geometric properties and related hypergeometric functions can be found in \cite{Carlson, Ponnu, Silverman}.

For analytic functions $f$ and $g$ in $\Delta$, we say that $f$ is subordinate to $g$, written as $f(z) \prec g(z), \; z\in \Delta,$ if there exists a Schwarz function $\omega$ (analytic in $\Delta$, $\omega(0)=0$, $|\omega(z)|<1$) such that $f(z)=g(\omega(z)).$
If $g$ is univalent, then $f(z) \prec g(z), \; z\in \Delta$ is equivalent to $f(0)=g(0)$ and $f(\Delta)\subset g(\Delta).$

Let $\mathcal{A}$ denote the class of analytic functions $f$ in $\Delta$  normalized by $f(0)=0$ and $f'(0)=1$. The subclass of univalent functions in $\mathcal{A}$ is denoted by $\mathcal{S}$. A function $f \in \mathcal{A}$ is called starlike if $\text{Re}\left(zf'(z)/f(z)\right)>0, \;z\in \Delta.$ The class of such functions is denoted by $\mathcal{S}^*$. Similarly, $f\in\mathcal{A}$ is called convex if $\text{Re}\left(1+zf''(z)/f'(z)\right)>0, \;z \in \Delta$. We denote the class of convex functions by $\mathcal{K}$. By Alexander's theorem \,$f\in\mathcal{K} \;\Leftrightarrow\; zf'(z)\in\mathcal{S}^{*}.$ A function $f\in\mathcal{A}$ is called close-to-convex if there exist a convex function $g\in\mathcal{A}$ and a real number $\theta$ such that
$\text{Re}\left(e^{i\theta} f'(z)/g'(z)\right)>0$. The class of close-to-convex functions is denoted by $\mathcal{C}$. Every close-to-convex function is univalent, but the converse is not generally true. Further details concerning these classes can be found in \cite{Duren,Goodman}.

For $0<\beta\leq1$, a function $f\in\mathcal{A}$ is said to be strongly starlike of order $\beta$ if
$$\bigg| \arg \bigg( \dfrac{zf'(z)}{f(z)} \bigg) \bigg| < \dfrac{\pi}{2}\beta, \qquad z \in \Delta,$$
Equivalently,
$$\frac{zf'(z)}{f(z)}
\prec\left(\frac{1+z}{1-z}\right)^\beta,
\qquad z\in\Delta,$$
where the principal branch is understood. The class of strongly starlike functions of order $\beta$, denoted by $\mathcal{SS}^*(\beta)$, was introduced independently by Stankiewicz \cite{Stan} and Brannan and Kirvan \cite{Brannan}. In particular, $\mathcal{SS}^*(1)=\mathcal{S}^*.$ A function $f\in\mathcal{A}$ is said to be strongly convex of order $\beta$, where $0<\beta\leq1$, if $zf'(z)$ is strongly starlike of order $\beta$. Equivalently,
$$\left|\arg\left(1+\frac{zf''(z)}{f'(z)}\right)\right|
<\frac{\pi\beta}{2},\qquad z\in\Delta.
$$
The corresponding class is denoted by $\mathcal{SK}(\beta)$. In particular, $\,\mathcal{SK}(1)=\mathcal{K}.$

The main objective of this article is to establish sufficient conditions for several geometric properties of confluent hypergeometric functions. We first establish lower bounds for the real part of $\Phi(a;c;z)$, subordination of its derivative, and strong convexity and strong starlikeness of prescribed order. We then investigate the close-to-convexity condition for $\Phi(a;c;z^m)$ and starlikeness of $z\Phi(a;c;z)$ through monotonicity properties of its coefficients. The results are compared with earlier results, and applications to related special functions are developed.

\medskip
We close this section with lemmas that will be used to prove several theorems.

\medskip
\noindent
{\bf Lemma 1.1} (\cite{miller2}).\;Let $\Lambda\subset\mathbb{C}$ and let a function $\Psi:\mathbb{C}^3\times \Delta\rightarrow \mathbb{C}$ satisfy the condition $\Psi(is,t,u+iv;z)\notin \Lambda$, for $z\in \Delta$ and for real $s$, $t$, $u$, $v$ satisfying $t\leq -(1+s^2)/2$ and $t+u\leq 0$. If $p$ is analytic in $\Delta$, with $p(0)=1$ and $\Psi(p(z),zp'(z),z^2p''(z);z)\in \Lambda$, $z\in \Delta$, then $\text{Re}(p(z))>0$ in $\Delta$.

\medskip
It may be noted that if we only have $\Psi:\mathbb{C}^2\times \Delta\rightarrow \mathbb{C}$, then the condition $\Psi(is,t;z)\notin \Lambda$ of Lemma 1.1 simplifies to $\Psi(is,t;z)\notin \Lambda$, for $z\in \Delta$ and for real $s$, $t$ with $t\leq -(1+s^2)/2$.

\medskip
\noindent
{\bf Lemma 1.2} (\cite{Fejer}). Let $f(z)=z+\sum_{n=2}^{\infty}A_nz^n$, with $A_1=1$. If $A_n\geq 0$, $\{nA_n\}_{n \ge 1}$ and $\{nA_n-(n+1)A_{n+1}\}_{n \ge 1}$ both are non-increasing, then $f \in \mathcal{S}^*$. 

\section{Main Results}

We first establish the following extension of \cite[Theorem 1.9]{Ponnu}.

\medskip
\noindent
{\bf Theorem 2.1} 
If $a, \,c, \,\beta$ are real such that $a\neq 0$ and 
$$c\ge
\begin{cases}
\displaystyle 1+\sqrt{1+\frac{(a-\beta(a+1))^2}{1-2\beta}}, & 0\le \beta <\dfrac{1}{2},\\
\displaystyle 1+\frac{|a|\beta}{1-\beta}
+\left|1-\frac{a\beta}{1-\beta}\right|, & \beta\ge\dfrac{1}{2}, \, \beta \neq 1.
\end{cases} \qquad \qquad \qquad $$
Then $\mbox{Re}\,\Phi(a; c; z)>\beta$ for all $ z \in \Delta.$
\begin{proof}
Let $p:\Delta\rightarrow \mathbb{C}$ be given by
\begin{equation}\label{norm.sub}
p(z)=\frac{\Phi(a;c;z)-\beta}{1-\beta}, \qquad (\beta \neq 1, \; z \in \Delta).
\end{equation}
Then $p$ is analytic in $\Delta$ and $p(0)=1$. As $\Phi(a;c;z)$ satisfies the differential equation \eqref{chf-2}, we get
\begin{equation}\label{CHF-starlike}
z^2 p''(z)+(c-z)zp'(z)-az \left(p(z)+\frac{\beta}{1-\beta} \right)=0.
\end{equation}
Let
$$ \Psi(w_1,w_2,w_3;z)=w_3+(c-z)w_2-az \left(w_1+\frac{\beta}{1-\beta} \right), $$
then \eqref{CHF-starlike} can be written as $\Psi(p(z),zp'(z),z^2p''(z);z)\in \Lambda$, where $\Lambda=\{ 0 \}$. By Lemma 1.1, $\text{Re} \,p(z)>0$ if $\Psi(is,t,u+iv;z) \notin \Lambda,$ where $s,\, t,\, u,\, v \in \mathbb{R}, \,t\leq-(1+s^2)/2,$ and $ t+u\leq0$  for $z \in \Delta$. We prove a stronger statement that $\text{Re} \,\Psi(i s,t,u+iv;z) <0$,  which eventually leads to the required condition. Considering $s,\, t,\, u,\, v \in \mathbb{R}$, $z=x+iy \in \Delta$, $t\leq-(1+s^2)/2$ and $ t+u\leq0$, we obtain
\begin{eqnarray*}
\text{Re}(\Psi(is,t,u+iv;z))
&=& u+t+(c-x-1)t+ays-\dfrac{ax\beta}{1-\beta} \\
&\leq & - \dfrac{1}{2} \left[(c-x-1)s^2-2ays+\left(c-x-1 +\dfrac{2ax\beta}{1-\beta}\right)\right] \notag\\
&:=&  Q_s(x,y,\beta). \notag
\end{eqnarray*}
Since the coefficient of $s^2$ in $Q_s(x,y,\beta)$ is negative, $Q_s(x,y,\beta)<0$ for all $s \in \mathbb{R}$ if its discriminant is strictly negative, that is,
$$(2ay)^2 < 4(c-x-1)\left(c-x-1 +\dfrac{2ax\beta}{1-\beta}\right).$$
Since $z \in \Delta,$ and so it is sufficient to verify the inequality $R_x(a,c,\beta)>0,$ where
$$R_x(a,c,\beta)=\left(a^2+1-\frac{2a\beta}{1-\beta}\right)x^2
-2(c-1)\left(1-\frac{a\beta}{1-\beta}\right)x
+(c-1)^2-a^2.$$

The case \(0\leq\beta<1/2\) was established in \cite[Theorem 1.9]{Ponnu}. We now consider the case for $\beta=1/2.$ If $\beta=1/2,$ then
\begin{eqnarray*}
R_x(a,c,1/2) &=& (a-1)^2x^2-2(c-1)(1-a)x+(c-1)^2-a^2\\
&=&\bigl((a-1)x+(c-1)\bigr)^2-a^2 \geq 0
\end{eqnarray*}
if and only if \,$c\ge 1+|a|+|a-1|\,$ for $a\neq0, \;c\geq1.$

Moreover, for $\beta \in (1/2, 1)$, we  set $t=\beta/(1-\beta), \; B=1-at,$ then $|t|>1$ and 
$$R_x(a,c,\beta)=\bigl(c-1-Bx\bigr)^2-a^2\bigl(1+(t^2-1)x^2\bigr).$$
Therefore,
$$R_x(a,c,\beta)\ge 0 \iff \bigl|c-1-Bx\bigr|\ge |a|\sqrt{1+(t^2-1)x^2}. $$

Hence $R_x(a,c,\beta)\ge 0$ for all $x\in(-1,1)$ if and only if either
$$c-1 \ge \sup_{|x|<1}\left(Bx+|a|\sqrt{1+(t^2-1)x^2}\right)$$
or
$$c-1 \le \inf_{|x|<1}\left(Bx-|a|\sqrt{1+(t^2-1)x^2}\right).$$

Now define $f(x)=Bx+|a|\sqrt{1+(t^2-1)x^2}.$ Since
$$f''(x)=|a|\frac{t^2-1}{\bigl(1+(t^2-1)x^2\bigr)^{3/2}}>0,$$
the function $f$ is convex on $(-1,1)$. Hence \,$\sup_{|x|<1} f(x)=\max\{f(1),f(-1)\}.$ But $f(1)=B+|a|t$ and $f(-1)=-B+|a|t,$ hence
$$\sup_{|x|<1} f(x) := \sup_{|x|<1} \left(Bx+|a|\sqrt{1+(t^2-1)x^2}\right) =|a|t+|B|.$$
Similarly,
$$\inf_{|x|<1}\left(Bx-|a|\sqrt{1+(t^2-1)x^2}\right)=-\bigl(|B|+|a|t\bigr).$$
Thus, the sufficient condition for $R_x(a,c,\beta)\ge 0$,
$$|c-1|\ge |a|t+|B|=\frac{|a|\beta}{1-\beta}+\left|1-\frac{a\beta}{1-\beta}\right|.$$
Since $c>1$, the required sufficient condition is
$$c\ge 1+\frac{|a|\beta}{1-\beta}
+\left|1-\frac{a\beta}{1-\beta}\right|.$$
This completes the proof.
\end{proof}

Now, using the first derivative identity of the confluent hypergeometric function
\begin{equation}\label{con-identity}
\Phi'(a;c;z)
=\frac{a}{c}\Phi(a+1;c+1;z)    
\end{equation}
and restricting the parameter $\beta\,$ to $\beta \in [0, 1)$, we obtain the following extension to the result \cite[Theorem 1.14]{Ponnu}. 

\medskip
\noindent
{\bf Corollary 2.2}
If $a, \,c, \,\beta$ are real such that $a\neq -1, \,\beta \in [0,1),$ and 
$$c\ge
\begin{cases}
\displaystyle \sqrt{1+\frac{(a+1-\beta(a+2))^2}{1-2\beta}}, & 0\le \beta <\dfrac{1}{2},\\
\displaystyle \frac{|a+1|\beta}{1-\beta}
+\left|1-\frac{(a+1)\beta}{1-\beta}\right|, & \dfrac{1}{2}\le\beta <1.
\end{cases}$$
Then the function $\dfrac{c}{a}(\Phi(a; c; z)-1)$ is univalent in $\Delta.$

\medskip
The generalized Laguerre polynomials $L_n^{(\alpha)}(z)$ are related with the confluent Hypergeometric function by the relation
$$\Phi(-n;\alpha+1;z) := M_{n, \alpha}(z)= \frac{n!}{(\alpha+1)_n}L_n^{(\alpha)}(z) \qquad (\alpha>-1, n \in \mathbb{N} \cup \{0\}).$$
Hence, taking $a=-n$ and $c=\alpha+1$ in Theorem 2.1, we obtain the following improved result of the normalized Laguerre polynomial $M_{n, \alpha}(z)$ \cite[Theorem 1]{praj1}.

\medskip
\noindent
{\bf Corollary 2.3}
Let $\alpha,\,\beta \in \mathbb{R}, \,n \in \mathbb{N}$. If
$$\alpha \ge
\begin{cases}
\displaystyle \sqrt{1+\frac{(n-\beta(n-1))^2}{1-2\beta}}, & 0\le \beta <\dfrac{1}{2},\\
\displaystyle 1+\frac{2n\beta}{1-\beta}, & \beta\ge\dfrac{1}{2}, \, \beta \neq 1.
\end{cases}$$
Then $\mbox{Re}\left(M_{n, \alpha}(z)\right)>\beta$ for all $ z \in \Delta.$

\medskip
Our next result gives conditions on $a$ and $c$ for which $(c/a) \Phi'(a;c;z) \prec ((1+z)/(1-z))^{1/2},$ or equivalently, the function $(c/a) \Phi'(a;c;z)$ lies in the sector $|\arg z|<\pi \beta/4, \;\beta \in (0, 1].$

\medskip
\noindent
{\bf Theorem 2.4}
Let $a, \,c, \beta \in \mathbb{R}, \,a \neq0,\,\beta \in (0,1].$ If 
$$c\ge \sqrt{1+\left(\frac{a+1}{\beta}\right)^2},\,
\quad \beta>0,$$
then $\dfrac{c}{a}\,\Phi'(a;c;z) \prec \left(\dfrac{1+z}{1-z}\right)^{\beta}.$

\begin{proof}
Let
$$p(z)=\exp\left[\frac1\beta\text{Log} \left(\dfrac{c}{a}\,\Phi'(a;c;z) \right)\right],$$
where $\text{Log}$ denotes the analytic logarithm in the right half-plane. Then $p$ is analytic in $\Delta$, $p(0)=1$, and
$$\dfrac{c}{a} \Phi'(a;c;z)=p(z)^\beta.$$
Substituting this identity into the differential equation \eqref{chf-2}, and multiplying by $z$, yields
\begin{eqnarray}\label{th3-1}
\qquad  \beta \,p(z) z^2p''(z)+\beta (\beta-1)(z \,p'(z))^2+\beta(c+1-z)p(z)(z\,p'(z))
-(a+1)z\,p(z)^2=0.  
\end{eqnarray}
This equation \eqref{th3-1} can be written as $\Psi(p(z), \,zp'(z), \,z^2p''(z);z)\in \Lambda$, where $\Lambda=\{0\}$. By Lemma 1.1, it is sufficient to show that $\Psi(is,t, u+iv;z) \notin \Lambda,$ where $s,t, u, v\in \mathbb{R}, \;t+u \leq 0$ and $t\leq-(1+s^2)/2$ for $z=x+iy\in \Delta$. Considering $s,t, u, v\in \mathbb{R}, \;z \in \Delta, \,t+u \leq 0$ and $t\leq-(1+s^2)/2$, we obtain
\begin{equation}\label{th3-2}
\text{Re} \,\Psi(is,t,u+iv;x+iy)=-\beta sv+\beta(\beta-1)t^2+\beta sty+(a+1)xs^2,
\end{equation}
and
\begin{equation}\label{th3-3}
\text{Im} \,\Psi(is,t,u+iv;x+iy)=\beta su+\beta st(c+1-x)+(a+1)ys^2.
\end{equation}

We now consider three cases depending upon $s=0, \,s>0, \,s<0.$

\medskip
\noindent
{\bf Case (i)}. If $s=0$, then from the equation \eqref{th3-2}, we have $\text{Re} \,\Psi(is,t,u+iv;z)=-\beta(1-\beta)t^2<0,$ thus $\text{Re} \,\Psi(is,t,u+iv;z) \notin \Lambda$.

\medskip
\noindent
{\bf Case (ii)}. If $s>0$, then using $u+t \leq 0$ and $t \leq -(1+s^2)/2,$ we get from the equation \eqref{th3-3},
\begin{equation}\label{th3-4}
\text{Im} \,\Psi(is,t,u+iv;x+iy)=s \left(-\frac{\beta}{2}(c-x)s^2+(a+1)ys-\frac{\beta}{2}(c-x) \right).    
\end{equation}
Let 
$$S_s(x,y,\beta)=-\frac{\beta}{2}(c-x)s^2+(a+1)ys-\frac{\beta}{2}(c-x)
\qquad (c>1,\,|z|<1,\,\beta\in[0,1)).$$
By hypothesis $\,c-x>0$, hence $S_s(x,y,\beta)$ is a quadratic polynomial in $s$ with negative leading coefficient. Therefore, $S_s(x,y,\beta)<0$ for all $\,s\in\mathbb{R}$ if and only if its discriminant is negative. Hence, the discriminant $D$ of $S_s(x,y,\beta)$,
\begin{eqnarray*}
D &=& (a+1)^2y^2-\beta^2(c-x)^2\\
&<&(a+1)^2(1-x^2)-\beta^2(c-x)^2.   
\end{eqnarray*}
Thus we need
$$\frac{1-x^2}{(c-x)^2}\le \frac{\beta^2
}{(a+1)^2}
\qquad \text{for all } -1<x<1.$$
Consider a function
$$f(x)=\frac{1-x^2}{(c-x)^2},\qquad -1<x<1,$$
which attains its maximum at $x=1/c$. Thus $\max_{-1<x<1} f(x)=f\!\left(1/c\right)=1/(c^2-1).$ Hence a necessary and sufficient condition for $\text{Im} \,\Psi(is,t,u+iv;x+iy)<0\,$ is $\;\dfrac{1}{c^2-1}\le \dfrac{\beta^2}{(a+1)^2},$ that is 
$$c\ge \sqrt{1+\left(\frac{a+1}{\beta}\right)^2},\,
\qquad \beta>0.$$

\medskip
\noindent
{\bf Case (iii)}. If $s<0$, then by using \eqref{th3-4}, with conditions $t\leq -(1+s^2)/2$ and $u+t \leq 0$ gives
$$\text{Im} \,\Psi(is,t,u+iv;x+iy)=s \left(-\frac{\beta}{2}(c-x)s^2+(a+1)ys-\frac{\beta}{2}(c-x) \right).$$
Since $s<0, \;\text{Im}\,\Psi >0$ provided $S_s(x,y,\beta)<0$, which reduces to case {\it(ii)} only.

Thus, in all the cases, $\Psi(is,t,u+iv;x+iy) \notin \Lambda$. This implies that
$$\frac{c}{a} \Phi'(a;c;z) \prec \left(\frac{1+z}{1-z}\right)^{\beta}, \qquad z \in \Delta.$$
This completes the proof.
\end{proof}

The next theorem establishes sufficient conditions for strong convexity of order $\beta$ for the confluent hypergeometric functions $\Phi(a;c;z)$.

\medskip
\noindent
{\bf Theorem 2.5}
Let $a, \,c, \beta \in \mathbb{R}, \,a \neq0$. Suppose that
$$\quad \alpha=\sin\left(\dfrac{\pi \beta}{2}\right), \; \gamma=\cos\left(\frac{\pi \beta}{2}\right), \; \beta_0=\dfrac{4}{\pi}\tan^{-1}\left(\dfrac{1}{2}\right), \; \eta_1(a, \beta)=1+\frac{|a|}{\alpha}
+\frac{1}{4\alpha^2}, \qquad \qquad \qquad \quad \quad $$
$$\eta_2(a, \beta)=1+\alpha^2 a^2
+\frac{\gamma |a|}{\alpha}
+\frac{1}{4\alpha^2}, \quad \text{and} \quad \eta_3(a, \beta)=1+\left(1+\frac{\gamma}{\alpha}\right)|a|. \qquad \qquad$$
If $\beta\in\left(0,\beta_0\right]$ and 
$$c\ge
\begin{cases}
\displaystyle
\eta_2(a, \beta), \;\;
& \text{if } |a|\le \dfrac{1}{\alpha(1+\gamma)},\\
\displaystyle
\eta_1(a, \beta),
& \text{if } \dfrac{1}{\alpha(1+\gamma)}\le |a| \le\dfrac{1}{2\alpha^2}\\
\displaystyle
\eta_3(a, \beta),
& \text{if } |a| >\dfrac{1}{2\alpha^2}.
\end{cases}$$
then $\Phi(a;c;z)$ is strongly convex of order $\beta.$ Moreover if  $\beta\in\left(\beta_0, 1\right)$ and  
$$c\geq
\begin{cases}
\eta_2(a,\beta),&0<|a|\leq \dfrac{1}{2\alpha^2},\\
\eta_3(a,\beta),&|a|>\dfrac{1}{2\alpha^2}.
\end{cases} \qquad \qquad $$
then $\Phi(a;c;z)$ is strongly convex of order $\beta.$

\begin{proof}
Let
$$q(z)=\exp\left[\frac1\beta\text{Log} \left(1+\frac{z \Phi''(z)}{\Phi'(z)} \right)\right],$$
where $\text{Log}$ denotes the analytic logarithm in the right half-plane. Then $q$ is analytic in $\Delta$, $q(0)=1$, and
\begin{eqnarray}\label{th2-2}
1+\frac{z\Phi''(z)}{\Phi'(z)}=q^\beta(z)   
\end{eqnarray}

Assume that $z \neq 0$. Dividing \eqref{chf-2} by $\Phi'(z)$ and using \eqref{th2-2}, we get
\begin{eqnarray}\label{th2-3}
a\frac{\Phi(z)}{\Phi'(z)}=q^\beta(z)+c-1-z.
\end{eqnarray}
Differentiate \eqref{th2-3} with respect to $z$, we get
\begin{eqnarray}\label{th2-35}
1-\frac{\Phi(z)}{\Phi'(z)}\cdot \frac{\Phi''(z)}{\Phi'(z)}=\frac{\beta q^{\beta-1}(z)q'(z)-1}{a}.
\end{eqnarray}
Further, using \eqref{th2-2} and \eqref{th2-3} in \eqref{th2-35}, we obtain 
\begin{eqnarray}\label{th2-4}
z \,q'(z)+\frac{1}{\beta}q^{\beta+1}(z)+\frac{c-z-2}{\beta}q(z)+\frac{1-c-az}{\beta}q^{1-\beta}(z)=0.
\end{eqnarray}
This equation \eqref{th2-4} can be written as $\Psi(w_1,w_2;z)\in \Lambda$, where $\Lambda=\{ 0 \}$. By Lemma 1.1, $\text{Re} \,p(z)>0$ if $\Psi(is,t;z) \notin \Lambda,$ where $s,\, t\in \mathbb{R}$ and $t\leq-(1+s^2)/2,$ for $z=x+iy\in \Delta$. We prove the stronger statement $\text{Re} \,\Psi(i s,t;z) <0$, which eventually leads to the required condition. Considering $s,\, t\in \mathbb{R}, \;z \in \Delta$ and $t\leq-(1+s^2)/2$, we obtain
\begin{eqnarray}\label{th2-5}
\Psi(is,t;z) =  t+\frac{1}{\beta}(is)^{\beta+1}+\frac{c-z-2}{\beta}(is)+\frac{1-c-az}{\beta}(is)^{1-\beta}.
\end{eqnarray}

We consider three cases depending upon $s=0,\,s>0$ and $s<0$. 

\medskip
\noindent
{\bf Case (i)}. If $s=0$, then $\text{Re} \,\Psi(is,t;z)=t<-(1+s^2)/2<0.$

\medskip
\noindent
{\bf Case (ii)}. If $s>0$, then using the principal branch,
$$(is)^{\beta+1}=
s^{\beta+1}e^{i(\beta+1)\pi/2},
\qquad
(is)^{1-\beta}=
s^{\,1-\beta}e^{i(1-\beta)\pi/2}.$$
For $s>0$ gives $\sqrt{is}=(1+i)\sqrt{s/2},$
which is the principal square root of $is$. Hence, the above power formulas are consistent with the principal branch. Now 
$$\text{Re}\!\left(\frac{1}{\beta}(is)^{\beta+1}\right)=-\frac{s^{\beta+1}}{\beta}\sin\frac{\beta\pi}{2}, \qquad \text{Re}\!\left(\frac{c-z-2}{\beta}(is)\right)
=\frac{ys}{\beta},$$
and
\begin{eqnarray*}
\text{Re}\!\left(\frac{1-c-az}{\beta}(is)^{1-\beta}\right)
&=&\frac{s^{1-\beta}}{\beta}
\left[(1-c-ax)\sin\frac{\beta\pi}{2}
+ay\cos\frac{\beta\pi}{2}\right].
\end{eqnarray*}
Therefore
\begin{eqnarray*}
\text{Re} \,\Psi(is,t;z)
&=&t-\frac{s^{\beta+1}}{\beta}\sin\frac{\beta\pi}{2}
+\frac{ys}{\beta}+\frac{s^{1-\beta}}{\beta}
\left[(1-c-ax)\sin\frac{\beta\pi}{2}
+ay\cos\frac{\beta\pi}{2}
\right]\\
& \leq &-\frac{1+s^2}{2} -\frac{s^{\beta+1}}{\beta}\sin\frac{\beta\pi}{2}
+\frac{ys}{\beta}+\frac{s^{1-\beta}}{\beta}
\left[(1-c-ax)\sin\frac{\beta\pi}{2}
+ay\cos\frac{\beta\pi}{2}
\right]\\
&:=& Q(s).
\end{eqnarray*}
Now we need to show  $Q(s)<0$ for all $s>0$. Let  $\alpha=\sin(\beta\pi/2)$ and $\gamma=\cos(\beta\pi/2).$ Then $\alpha>0$, and
\begin{equation}\label{Q(s)}
Q(s)= -\frac{1+s^2}{2} -\frac{\alpha}{\beta}s^{\beta+1}
+\frac{ys}{\beta}+\frac{s^{1-\beta}}{\beta}
\Bigl[(1-c-ax)\alpha+ay\gamma\Bigr]
\end{equation}

Since $|x|,|y|<1$, we have $ys< s,$ and by the Cauchy–Schwarz inequality, we have
$$-ax\alpha+ay\gamma \le |a|\sqrt{\alpha^2+\gamma^2}\sqrt{x^2+y^2}
< |a|.$$
Also
$$s^\beta-\alpha \,s^{2\beta}=-\alpha \left(s^\beta-\frac{1}{2\alpha} \right)^2+\frac{1}{4\alpha} < \frac{1}{4\alpha}.$$
Hence
\begin{eqnarray*}
Q(s) &<&-\frac{1+s^2}{2}-\frac{\alpha}{\beta}s^{\beta+1}+\frac{s}{\beta}
+\frac{s^{1-\beta}}{\beta}\bigl((1-c)\alpha+|a|\bigr)\\
&=& -\frac{1+s^2}{2}
+\frac{s^{1-\beta}}{\beta}
\Bigl[s^\beta-\alpha s^{2\beta}+(1-c)\alpha+|a|\Bigr]\\
&<&-\frac{1+s^2}{2}
+\frac{s^{1-\beta}}{\beta}
\left[(1-c)\alpha+|a|+\frac{1}{4\alpha}\right].\\
&\leq& 0, \qquad \mbox{if} \quad (1-c)\alpha+|a|+\frac{1}{4\alpha}\le 0.
\end{eqnarray*}
Therefore, a sufficient condition for $Q(s)<0$ is
\begin{equation}\label{first-c}
c\ge 1+\frac{|a|}{\alpha}
+\frac{1}{4\alpha^2}:=\eta_1(a, \beta)
\end{equation}

\medskip
Now using a different inequality, we obtain another lower bound for $c$. Using \eqref{Q(s)}, we have
\begin{eqnarray*}
Q(s) &=& -\frac{1+s^2}{2} + \frac{s^{1-\beta}}{\beta}
\Bigl(-\alpha s^{2 \beta}+y s^{\beta}+(1-c-ax)\alpha+ay\gamma\Bigr)\\
&:=& -\frac{1+s^2}{2} + \frac{s^{1-\beta}}{\beta} R(s^{\beta}),
\end{eqnarray*}
where $R(u)=-\alpha u^2+yu+(1-c-ax)\alpha+ay\gamma.$ Now the discriminant $D$ of $R(u)$
\begin{eqnarray*}
D &=& y^2+4 \alpha \left((1-c-ax)\alpha+ay\gamma \right)\\
&<& 1-x^2+4 \alpha^2(1-c-ax)+4 \alpha \gamma |a|:= -S(x), \qquad \qquad \qquad \qquad 
\end{eqnarray*}
where
\begin{eqnarray*}
S(x) &=& x^2+4\alpha^2ax+4\alpha^2c-4\alpha^2-4\alpha\gamma |a|-1\\
&=& \left(x+2\alpha^2 a\right)^2-4\alpha^4 a^2+4\alpha^2c-4\alpha^2-4\alpha\gamma |a|-1.
\end{eqnarray*}
Since $S(x)$ is a quadratic polynomial in $x$ with positive leading coefficient, its minimum over $(-1,1)$ is attained either at the vertex or at the boundary points $x=\pm 1$. The vertex of $S(x)$ is $\,x_0=-2\alpha^2 a.$ Thus, we consider two cases.

Let $\;|x_0|\le 1$, that is, $\,|a|\le 1/2\alpha^2.$ Then the vertex lies in $[-1,1]$, hence
$$\min S(x)=S(x_0)=-4\alpha^4 a^2+4\alpha^2c-4\alpha^2-4\alpha\gamma |a|-1.$$
Therefore $S(x)\ge 0$ for all $x\in(-1,1)$ if and only if
$$c\ge \alpha^2a^2+1+\frac{\gamma}{\alpha}|a|+\frac{1}{4\alpha^2}.$$

Moreover, let $\;|x_0|>1$, that is, $\,|a|>1/2\alpha^2.$ Then the vertex lies outside $(-1,1)$, hence
$$\inf_{x\in(-1,1)} S(x)=
1-4\alpha^2|a|+4\alpha^2c-4\alpha^2-4\alpha\gamma |a|-1.$$
Therefore, $S(x)\ge 0$ for all $x\in(-1,1)$ if and only if
$$c\ge 1+|a|+\frac{\gamma}{\alpha}|a|.$$
Combining these two cases, the best condition on $c$ is
$$c\ge
\begin{cases}
\displaystyle
1+\alpha^2 a^2
+\frac{\gamma |a|}{\alpha}
+\frac{1}{4\alpha^2}:=\eta_2(a, \beta), \;\;
& \text{if } |a|\le \dfrac{1}{2\alpha^2},\\
\displaystyle
1+\left(1+\frac{\gamma}{\alpha}\right)|a|:=\eta_3(a, \beta),
& \text{if } |a|> \dfrac{1}{2\alpha^2}.
\end{cases}$$

We now compare different lower bounds of $c$. If $\beta \in (0, \beta_0],$ then $\dfrac{1}{\alpha(1+\gamma)} \le\dfrac{1}{2\alpha^2}$. Using $1-\gamma=\alpha^2/(1+\gamma)$, we have
$$\eta_2-\eta_1=\alpha^2|a|^2+\frac{\gamma-1}{\alpha}|a|=|a|\alpha\left(\alpha |a|-\frac{1}{1+\gamma}\right).$$
Hence
$$\eta_2 \le \eta_1\iff |a| \le \frac{1}{\alpha(1+\gamma)} \quad \text{and} \quad \eta_2>\eta_1 \iff |a|>\frac{1}{\alpha(1+\gamma)}.$$
Therefore
$$c\ge
\begin{cases}
\displaystyle
\eta_2(a, \beta), \;\;
& \text{if } |a|\le \dfrac{1}{\alpha(1+\gamma)},\\
\displaystyle
\eta_1(a, \beta),
& \text{if } \dfrac{1}{\alpha(1+\gamma)}\le |a| \le\dfrac{1}{2\alpha^2}\\
\displaystyle
\eta_3(a, \beta),
& \text{if } |a| >\dfrac{1}{2\alpha^2}.
\end{cases}$$

Moreover, if $\beta \in (\beta_0, 1]$, then $\,\dfrac{1}{2\alpha^2}<\dfrac{1}{\alpha(1+\gamma)}.$ Hence
$$c\geq
\begin{cases}
\eta_2(a,\beta),&0<|a|\leq \dfrac{1}{2\alpha^2},\\
\eta_3(a,\beta),&|a|>\dfrac{1}{2\alpha^2}.
\end{cases} \qquad \qquad $$

\medskip
\noindent
{\bf Case (iii).} If $s<0.$ In this case $\sqrt{is}=(1-i)\sqrt{m/2}$, where $m=-s>0$. Again from \eqref{th2-5}, we get
\begin{eqnarray*}
\text{Re} \,\Psi(is,t;z)
&=&t-\frac{m^{\beta+1}}{\beta}\sin\frac{\beta\pi}{2}
+\frac{ym}{\beta}+\frac{m^{1-\beta}}{\beta}
\left[(1-c-ax)\sin\frac{\beta\pi}{2}
+ay\cos\frac{\beta\pi}{2}
\right]\\
& \leq &-\frac{1+m^2}{2} -\frac{m^{\beta+1}}{\beta}\sin\frac{\beta\pi}{2}
+\frac{ym}{\beta}+\frac{m^{1-\beta}}{\beta}
\left[(1-c-ax)\sin\frac{\beta\pi}{2}
+ay\cos\frac{\beta\pi}{2}
\right].
\end{eqnarray*}
The last inequality is the same as the one in case (ii) and therefore holds.

Combining all the cases and using Lemma 1.1, we conclude that $\Phi(a; c; z)$ is strongly convex of order $\beta$. 
\end{proof}

Using the derivative identity \eqref{con-identity}, we get
$$1+\frac{\Phi''(a-1; c-1; z)}{\Phi'(a-1; c-1; z)}=\frac{z \left(\Phi'(a; c; z) \right)'}{z \Phi(a; c; z)}.$$
Hence, we obtain the following result of strong starlikeness. 

\medskip
\noindent
{\bf Corollary 2.6}
Let $a, \,c, \beta \in \mathbb{R}, \,a \neq 1$. Suppose that
$$\quad \alpha=\sin\left(\dfrac{\pi \beta}{2}\right), \; \gamma=\cos\left(\frac{\pi \beta}{2}\right), \; \beta_0=\dfrac{4}{\pi}\tan^{-1}\left(\dfrac{1}{2}\right), \; \tilde{\eta}_1(a, \beta)=1+\frac{|a-1|}{\alpha}
+\frac{1}{4\alpha^2}, \qquad \qquad \qquad $$
$$\tilde{\eta_2}(a, \beta)=1+\alpha^2 (a-1)^2+\frac{\gamma |a-1|}{\alpha}
+\frac{1}{4\alpha^2}, \;\; \text{and} \;\; \tilde{\eta_3}(a, \beta)=1+\left(1+\frac{\gamma}{\alpha}\right)|a-1|. $$
If $\beta\in\left(0,\beta_0\right]$ and 
$$c\ge
\begin{cases}
\displaystyle
1+\tilde{\eta_2}(a, \beta), \;\;
& |a-1|\le \dfrac{1}{\alpha(1+\gamma)},\\
\displaystyle
1+\tilde{\eta_1}(a, \beta),
& \dfrac{1}{\alpha(1+\gamma)}\le |a-1| \le\dfrac{1}{2\alpha^2}\\
\displaystyle
1+\tilde{\eta_3}(a, \beta),
& |a-1| >\dfrac{1}{2\alpha^2}.
\end{cases}$$
then $\Phi(a;c;z)$ is strongly starlike of order $\beta.$ Moreover if  $\beta\in\left(\beta_0, 1\right)$ and  
$$c\geq
\begin{cases}
1+\tilde{\eta_2}(a,\beta),&0<|a-1|\leq \dfrac{1}{2\alpha^2},\\
1+\tilde{\eta_3(a,\beta)},&|a-1|>\dfrac{1}{2\alpha^2}.
\end{cases} \qquad \qquad $$
then $\Phi(a;c;z)$ is strongly starlike of order $\beta.$

\medskip
Setting $\beta=1/2$ in Theorem 2.5 and Corollary 2.6 yields the following improvement in \cite[Theorem 1.1, Corollary 1.2]{nbr}.   

\medskip
\noindent
{\bf Corollary 2.7}
Let $a,c\in\mathbb R$, $a\neq0$. Suppose that
$$c\geq
\begin{cases}
\displaystyle
\frac32+|a|+\frac{a^2}{2},
&0<|a|\leq2(\sqrt2-1),\\
\displaystyle
\frac32+\sqrt2|a|,
&2(\sqrt2-1)<|a|\leq1,\\
\displaystyle
1+2|a|, &|a|>1.
\end{cases}$$
Then $\Phi(a;c;z)$ is strongly convex of order $1/2$ in the unit disk $\Delta$.

\medskip
\noindent
{\bf Corollary 2.8}
Let $a,c\in\mathbb R$, $a\neq1$. Suppose that
$$c\geq
\begin{cases}
\displaystyle
\frac{5}{2}+|a-1|+\frac{(a-1)^2}{2},
&0<|a-1|\leq2(\sqrt2-1),\\
\displaystyle
\frac{5}{2}+\sqrt2|a-1|,
&2(\sqrt2-1)<|a-1|\leq1,\\
\displaystyle
2+2|a-1|, &|a-1|>1.
\end{cases}$$
Then $\Phi(a;c;z)$ is strongly starlike of order $1/2$ in the unit disk $\Delta$.

\medskip
To prove the next result, we need the following Lemma.

\medskip
\noindent
{\bf Lemma 2.9} 
Let $m \in \mathbb{N}$ and let $\,f(z)=z+\sum_{n=1}^{\infty}a_{mn+1}z^{mn+1}\in\mathcal A.$ Suppose that the coefficients are real and satisfy either
$$1\geq(m+1)a_{m+1}
\geq(2m+1)a_{2m+1}
\geq\cdots\geq0, \qquad \qquad \qquad \qquad \qquad \qquad $$
or
$$1\leq (m+1)a_{m+1}
\leq (2m+1)a_{2m+1}
\leq \cdots,\qquad \lim_{n\to\infty}(mn+1)a_{mn+1}\leq2.$$
Then $f$ is close-to-convex with respect to $g_m(z)=\int_0^z (1-t^m)^{-1}\,dt \,$ in $\,\Delta\,$ for $\,m\leq2$, and in $|z|<R_m:=(m-1)^{-1/m}$ for $m\geq3$.

\begin{proof}
Clearly $g_m \in \mathcal{A}$ and
\begin{equation*}
1+\frac{zg_m''(z)}{g_m'(z)}
=\frac{1+(m-1)z^m}{1-z^m}.
\end{equation*}
Thus $g_1$ and $g_2$ are convex in $\Delta$. Now suppose that $m\geq3$. Since $m>2$ and $\text{Re}(z^m)\geq-|z|^m$, then
$$\begin{aligned}
\text{Re}\left(
1+\frac{zg_m''(z)}{g_m'(z)}
\right)
&\geq\frac{
1-(m-2)|z|^m-(m-1)|z|^{2m}
}{|1-z^m|^2}\\
&=\frac{(1+|z|^m)(1-(m-1)|z|^m)
}{|1-z^m|^2}>0
\end{aligned}$$
whenever $\,|z|<(m-1)^{-1/m}=R_m.$ Hence $g_m$ is convex in $|z|<R_m$.

In either case (i) or (ii), the sequence $\{b_n\}$ is monotone and satisfies $\,\sum_{n=1}^{\infty}|b_n-b_{n-1}|
=\left|\lim_{n\to\infty}b_n-b_0\right|
\leq1.$ Hence, for every $z\in\Delta$,
$$\begin{aligned}
\text{Re}\left\{
\frac{f'(z)}{g_m'(z)}\right\} &=\text{Re}\left\{(1-z^m)f'(z)\right\}\\
&\geq 1-\sum_{n=1}^{\infty}
|b_n-b_{n-1}||z|^{mn}>0.
\end{aligned}$$
Therefore, $f$ is close-to-convex with respect to $g_m$ in the respective disks.
\end{proof}

We now prove a sufficient condition for the close-to-convexity of $z \,\Phi(a;c;z^m)$.

\medskip
\noindent
{\bf Theorem 2.10.}
Let $m\in \mathbb{N}, \;a>0$, and $c\geq(m+1)a$. Then the function $\,f_m(z)=z \,\Phi(a;c;z^m)$ is close-to-convex in $\Delta$ with respect to $\,g_m(z)=\int_0^z(1-t^m)^{-1} dt$ for $m \le 2,$ and close-to-convex in $|z|<R_m:=(m-1)^{-1/m}$ for $m \ge 3.$

\begin{proof}
Let 
$$f_m(z)=\ z \,\Phi(a;c;z^m)=z+\sum_{n=1}^{\infty}A_{n}z^{mn+1},$$ where $A_n=(a)_n/(c)_n \,n!$ and set $b_n=(mn+1)A_n$, where $A_0=1$. Then $b_0=1, \;b_1 \le 1,$ and
$$\frac{b_{n}}{b_{n+1}}-1=\frac{(mn+1)(c+n)(n+1)-(mn+m+1)(a+n)}
{(mn+1)(c+n)(n+1)}.$$
For $n\geq0$, the hypothesis $c\geq(m+1)a$
gives
\begin{eqnarray*}
D_n &=&(mn+1)(c+n)(n+1)-(mn+m+1)(a+n)\\
& \geq & a\left[(m+1)(mn+1)(n+1)-(mn+m+1)\right]+n\left[(mn+1)(n+1)-(mn+m+1)\right]\\
&=&a\left[m(m+1)n^2+\left(m^2+m+1\right)n\right] +n\left(mn^2+n-m \right) \ge 0.
\end{eqnarray*}

This implies that $b_n \ge b_{n+1}, \;n \ge 1.$ Hence $\,1=b_0\geq b_1\geq b_2\geq\cdots>0.$ Now using Lemma 2.9, $f_m$ is close-to-convex with respect to $g_m$ in respective discs. 
\end{proof}

\begin{remark}
(i) \,Taking $m=1$ in Theorem~2.10, we conclude that $z\Phi(a;c;z)$ is close-to-convex with respect to $-\log(1-z)$ in $\Delta$, provided that $a>0$ and $c \ge2a$. For $0<a<1/3$, this condition
is best than the corresponding sufficient condition $c\geq (1+a)/2$ obtained in \cite[Theorem~1.3]{Ponnu}.
\item[(ii)] Taking $m=2$ in Theorem~2.10, we conclude that $z\Phi(a;c;z^2)$ is close-to-convex with respect to $(1/2)\log[(1+z)/(1-z)]$ in $\Delta$, provided that $a>0$ and $c\ge3a$. For $\,0<a<2/7$, this condition
is best than the corresponding sufficient condition $c \geq 2(1+a)/3\,$ obtained in \cite[Theorem~1.5]{Ponnu}. 
\item[(iii)] Taking $m=3$ in Theorem~2.10, we conclude that $z\Phi(a;c;z^3)$ is close-to-convex with respect to 
$$g_3(z)= -\frac{1}{3}\log(1-z)
+\frac{1}{6}\log(1+z+z^2)
+\frac1{\sqrt3}
\left(\arctan\frac{2z+1}{\sqrt3}
-\frac{\pi}{6}\right)$$
in $|z|<2^{-1/3}\approx 0.7937$, provided that $a>0$ and $c\ge 4a$. 
\item[(iv)] Taking $m=4$ in Theorem~2.10, we conclude that $z\Phi(a;c;z^4)$ is close-to-convex with respect to 
$$g_4(z)=\frac{1}{4}\log\frac{1+z}{1-z}
+\frac{1}{2}\arctan z$$
in $|z|<3^{-1/4}\approx 0.7598$, provided that $a>0$ and $c \ge 5a$. 
\end{remark} 

\medskip
Finally, using Lemma 1.2 and Theorem 2.10, we prove the following improved version of \cite[Theorem 1.8]{Ponnu}.

\medskip
\noindent
{\bf Theorem 2.11}
Let $a>0$, and define
$$C_1(a)=
\frac{4a-1+\sqrt{10a^2+2a+1}}{2} \quad \text{and} \quad C_2(a)=
\frac{9a-9+\sqrt{3(11a^2+42a+43)}}{12}.$$
If $\,c\geq\max\{C_1(a),C_2(a)\},$ then $\,f(z)=z\Phi(a;c;z)$ is starlike and close-to-convex with respect to
$-\log(1-z)$ in $\Delta$.

\begin{proof}
Let $\,f(z)=z\Phi(a;c;z)
=z+\sum_{n=2}^{\infty}A_nz^n,$ where
$$A_1=1,\qquad
A_n=\frac{(a)_{n-1}}{(c)_{n-1}(n-1)!},
\quad n\geq2.$$
Since $a>0$ and $c\geq C_1(a)>2a>0$,
we have $A_n>0$ for every $n\geq1$. By Theorem~2.10, the sequence
$\{nA_n\}_{n\geq1}$ is non-increasing.
Hence, $f$ is close-to-convex
with respect to $-\log(1-z)$ in
$\Delta$.

To establish starlikeness, define $\,B_n=nA_n-(n+1)A_{n+1},
\, n\geq1.$ We shall show that $\{B_n\}_{n\geq1}$ is non-increasing. First note that
$$B_1-B_2=\frac{
2c^2+(2-8a)c+3a^2-5a
}{2c(c+1)} \ge 0$$
for all $c\geq C_1(a)$, where $C_1(a)$
is the larger root of the quadratic
polynomial in the numerator.

Next, for $n\geq2$, we obtain
$$B_n-B_{n+1}=\frac{P(n)}
{n(n+1)(c+n-1)(c+n)} A_n,$$
where $P$ is a polynomial of degree
five. Writing $n=y+2$, we obtain
$$P(y+2)=y^5+p_4\,y^4+p_3\,y^3+p_2\,y^2+p_1\,y+p_0,$$
with
$$\begin{aligned}
p_4&=2(c+4), \quad p_3=c^2+15c-2a+22, \quad p_2=7c^2+39c-2ac-14a+23,\\
p_1&=a^2+16c^2+42c-12ac-31a+4, \quad p_0=2(6c^2+9c-9ac+2a^2-12a-2).
\end{aligned}$$
We now verify that $p_j\ge0$ for $j=0,1,\ldots,4$.

Clearly, $p_4>0$. Moreover,
$p_1,p_2$, and $p_3$ are increasing
functions of $c$ for $c\geq2a$.
Evaluating them at $c=2a$, we find
$$p_3>p_3\big|_{c=2a}, \qquad p_1\ge p_1\big|_{c=2a}>0, \qquad \text{and} \qquad p_2\ge p_2\big|_{c=2a}>0.$$
Hence $p_1,p_2,p_3>0$ whenever $c\geq2a$.

It remains to consider $p_0$. Write $\,p_0=2Q(c),$ where $Q(c)=6c^2+(9-9a)c
+2a^2-12a-2.$ Since $Q$ is strictly increasing on $[2a,\infty)$ and $C_2(a)$ is the larger root of $Q(c)=0$, the condition $c\geq\max\{C_1(a),C_2(a)\}$
implies $\,Q(c)\geq0.$ Therefore $p_0\geq0$. It follows that $\,P(y+2)\geq0,\; y\geq0.$ In particular, $P(n)\geq0$ for
every integer $n\geq2$. This implies $\,B_n-B_{n+1}\ge0 \; (n\ge2).$
Now combining this inequality with $B_1\ge B_2$, we conclude that $\{B_n\}_{n\ge1}$ is non-increasing. Therefore Lemma~1.2 yields $f(z)=z\Phi(a;c;z)\in\mathcal S^*.$ This completes the proof.
\end{proof}

\medskip
\noindent
{\bf Concluding remarks.}
The close connections between confluent hypergeometric functions and several classical special functions, including generalized Laguerre polynomials, Whittaker functions, Bessel functions, and Hermite polynomials, provide opportunities to extend the results obtained in this paper. Through suitable parameter specializations, the established sufficient conditions may be used to investigate the geometric properties of these related special functions.

Furthermore, the present study opens avenues for future research, particularly in determining sharper parameter bounds and the exact radii of starlikeness, convexity, and close-to-convexity for the associated functions. Such investigations may lead to refined geometric criteria and a deeper understanding of the relationships between confluent hypergeometric functions and their special cases.

\medskip
\noindent
{\bf Funding Declaration.} The authors declare that no financial support was received for this manuscript.

\medskip
\noindent
{\bf Clinical Trial Declaration.} This manuscript does not involve any clinical trials.

\end{document}